
\catcode`\@=11
\def\undefine#1{\let#1\undefined}
\def\newsymbol#1#2#3#4#5{\let\next@\relax
 \ifnum#2=\@ne\let\next@\msafam@\else
 \ifnum#2=\tw@\let\next@\msbfam@\fi\fi
 \mathchardef#1="#3\next@#4#5}
\def\mathhexbox@#1#2#3{\relax
 \ifmmode\mathpalette{}{\m@th\mathchar"#1#2#3}%
 \else\leavevmode\hbox{$\m@th\mathchar"#1#2#3$}\fi}
\def\hexnumber@#1{\ifcase#1 0\or 1\or 2\or 3\or 4\or 5\or 6\or 7\or 8\or
 9\or A\or B\or C\or D\or E\or F\fi}

\newdimen\ex@
\ex@.2326ex
\def\varinjlim{\mathop{\vtop{\ialign{##\crcr
 \hfil\rm lim\hfil\crcr\noalign{\nointerlineskip}\rightarrowfill\crcr
 \noalign{\nointerlineskip\kern-\ex@}\crcr}}}}
\def\varprojlim{\mathop{\vtop{\ialign{##\crcr
 \hfil\rm lim\hfil\crcr\noalign{\nointerlineskip}\leftarrowfill\crcr
 \noalign{\nointerlineskip\kern-\ex@}\crcr}}}}
\def\varliminf{\mathop{\underline{\vrule height\z@ depth.2exwidth\z@
 \hbox{\rm lim}}}}

\font\tenmsa=msam10
\font\sevenmsa=msam7
\font\fivemsa=msam5
\newfam\msafam
\textfont\msafam=\tenmsa
\scriptfont\msafam=\sevenmsa
\scriptscriptfont\msafam=\fivemsa
\edef\msafam@{\hexnumber@\msafam}
\mathchardef\dabar@"0\msafam@39
\def\dashrightarrow{\mathrel{\dabar@\dabar@\mathchar"0\msafam@4B}}
\def\dashleftarrow{\mathrel{\mathchar"0\msafam@4C\dabar@\dabar@}}

\font\tenmsb=msbm10
\font\sevenmsb=msbm7
\font\fivemsb=msbm5
\newfam\msbfam
\textfont\msbfam=\tenmsb
\scriptfont\msbfam=\sevenmsb
\scriptscriptfont\msbfam=\fivemsb
\edef\msbfam@{\hexnumber@\msbfam}
\def\Bbb#1{{\fam\msbfam\relax#1}}
\def\widehat#1{\setbox\z@\hbox{$\m@th#1$}%
 \ifdim\wd\z@>\tw@ em\mathaccent"0\msbfam@5B{#1}%
 \else\mathaccent"0362{#1}\fi}
\font\teneufm=eufm10
\font\seveneufm=eufm7
\font\fiveeufm=eufm5
\newfam\eufmfam
\textfont\eufmfam=\teneufm
\scriptfont\eufmfam=\seveneufm
\scriptscriptfont\eufmfam=\fiveeufm

\newsymbol\boxtimes 1202

\catcode`\@=12

\magnification=\magstep1
\font\title = cmr10 scaled \magstep2

\font\smallmath= cmmi8
\font\smalltext = cmr7
\font\tinymath=cmmi5
\font\smallsym = cmsy7
\font\author = cmcsc10
\font\addr = cmti8
\font\byabs = cmr8

\font\email = cmtt8

\parindent=1em
\baselineskip 15pt
\vsize=18.5 cm

\newcount\refcount
\newcount\seccount
\newcount\sscount
\newcount\eqcount
\newcount\boxcount
\newcount\testcount
\newcount\bibcount
\boxcount = 128
\seccount = -1

\def\sec#1{\advance\seccount by 1\bigskip\goodbreak\noindent
	{\bf\number\seccount.\ #1}\medskip \sscount = 0\eqcount = 0}
\def\proc#1#2{\advance\sscount by 1\eqcount = 0
	\medskip\goodbreak\noindent{\author #1}
	{\tenrm{\number\sscount}}:\ \ {\it #2}}
\def\nproc#1#2#3{\advance\sscount by 1\eqcount = 0\global
	\edef#1{#2\ \number\sscount}	
	\medskip\goodbreak\noindent{\author #2}
	{\tenrm{\number\sscount}}:\ \ {\it #3}}
\def\proof{\medskip\noindent{\it Proof:\ \ }}

\def\eql#1{\global\advance\eqcount by 1\global
	\edef#1{(\number\sscount.\number\eqcount)}\leqno{#1}}
\def\ref#1#2{\advance\refcount by 1\global
	\edef#1{[\number\refcount]}\setbox\boxcount=
	\vbox{\item{[\number\refcount]}#2}\advance\boxcount by 1}
\def\biblio{{\frenchspacing
	\bigskip\goodbreak\centerline{\it References}\medskip
	\bibcount = 128\loop\ifnum\testcount < \refcount
	\goodbreak\advance\testcount by 1\box\bibcount
	\advance\bibcount by 1\vskip 4pt\repeat\medskip}}
\def\tightmatrix#1{\null\,\vcenter{\normalbaselines\mathsurround=0pt
	\ialign{\hfil$##$\hfil&&\ \hfil$##$\hfil\crcr
	\mathstrut\crcr\noalign{\kern-\baselineskip}
	#1\crcr\mathstrut\crcr\noalign{\kern-\baselineskip}}}\,}

\ref\CN{I.~L.~Chuang, M.~A.~Nielsen: {\it Quantum computation and quantum information.}   
    Cambridge University Press, Cambridge, 2000.}
\ref\HW{G.~H.~Hardy, E.~M.~Wright:  {\it An Introduction to the Theory of 
    Numbers.} fourth edition, Oxford University Press, Oxford, 1960.}
\ref\K{D.~E.~Knuth: {\it The Art of Computer Programming. Vol. 2: Seminumerical algorithms.}
    Addison-Wesley Publishing Co., Reading, MA, 1969.}
\ref\KY{D.~E.~Knuth, A.~C.~Yao:
    Analysis of the subtractive algorithm for greatest common divisors.  
    {\it Proc.\ Nat.\ Acad.\ Sci.\ U.S.A.} {\bf 72}  (1975), no. 12, 4720--4722.}
\ref\Lind{S.~Lindhurst: An analysis of Shanks's algorithm 
    for computing square roots in finite fields.
    {\it Number theory (Ottawa, ON, 1996)},  231--242,
    CRM Proc.\ Lecture Notes, 19,
    Amer.\ Math.\ Soc., Providence, RI, 1999.}
\ref\Lu{A.~Lubotzky: {\it Discrete groups, expanding graphs and invariant 
    measures.} Birkh\"auser Verlag, Basel, 1994.}
\ref\LPS{A.~Lubotzky, R.~Phillips, P.~Sarnak: Ramanujan graphs,
    {\it Combinatorica} {\bf 8} (1988) 261--277.}
\ref\Shanks{A.~Tonelli: Sulla risoluzione della congruenza 
    $x^2 \equiv c$\ (mod $p^\lambda$), 
    {\it Atti R.\ Accad.\ Lincei} {\bf 1} (1892), 116--120.}

\def\emph{\bf}

\def\|{|\;}


\def\Z{{\Bbb Z}}

\def\F{{\Bbb F}}

\def\N{{\Bbb N}}

\def\SL{{\rm SL}}

\def\qed{\hfill\hbox{$\sqcup$\llap{$\sqcap$}}\medskip}


\centerline{\bf Navigating the
Cayley graph of $\SL_2 (\F_p)$}
\bigskip
\centerline{\byabs by}
\medskip
\noindent{\author\hfill Michael Larsen\footnote*
{\sevenrm Partially supported by NSF
Grant DMS 97-27553.}\hfill}
\medskip
\centerline{\addr Department of Mathematics, Indiana University}
\centerline{\addr Bloomington, IN 47405, USA\footnote{}{AMS Classification
20F05}}
\centerline{\email larsen@math.indiana.edu}
\medskip
\centerline{\byabs ABSTRACT}
\smallskip
{\byabs \narrower\narrower
\textfont0 = \smalltext
\textfont1 = \smallmath
\scriptfont1 = \tinymath
\textfont2 = \smallsym
We present a non-deterministic polynomial-time algorithm to find a path of length
$O(\log p \log\log p)$ between any two vertices of the Cayley graph of 
$\SL_{\hbox{\tinymath 2}}(\hbox{{\sevenmsb F}}_p)$.}
\medskip

\def\Z{{\Bbb Z}}

\def\N{{\Bbb N}}

\def\qed{\hfill\nobreak\rlap{$\sqcup$}$\sqcap$}

\def\F{{\Bbb F}}

\def\mx#1#2#3#4{\left(\matrix{#1&#2\cr#3&#4\cr}\right)}

It is well known that $\SL_2(\F_p)$ is generated by $\mx1101$ and 
$\mx1011$.  It is a much deeper theorem \Lu\ 
that the Cayley diameter of this group 
with respect to these generators is $O(\log p)$.  
There are two known proofs.  One
depends on uniformly bounding the eigenvalues of the Laplacian on 
$L^2_0(X(p))$ away from zero \Lu. 
The other uses the circle method to show that any element of $\SL_2(\F_p)$ lifts
to an element of $\SL_2(\Z)$ which has a short word representation \LPS.  
Neither method is constructive. 
A.~Lubotzky asked \Lu\  for an 
efficient algorithm to find short word representations of general elements of
$\SL_2(\F_p)$.
In this note we give such an 
algorithm, but for word representations
of length $O(\log p \log\log p)$ rather than $O(\log p)$.
More precisely, we prove 
\proc{Theorem}{There exist constants $c_1$ and $c_2$ such that
for any $c_3 < 1$, there exists $c_4$
such that for any prime $p$ and any element of $\SL_2(\F_p)$,
the algorithm will find a word of length $\le c_1 \log p \log\log p$ 
in time $\le c_4 \log^{c_2} p$
with probability $\ge c_3$.}
\smallskip
 
\medskip
Consider first the basic strategy of lifting
$\alpha\in\SL_2(\F_p)$ to 
$\tilde\alpha\in\SL_2(\Z)$ 
and then using Euclid's algorithm to represent $\tilde\alpha$.
The trouble is that we must use the {\emph subtractive} Euclidean algorithm.
That is, we have to pay for each operation of 
subtracting
one row from another, so the performance of the algorithm is worse 
than that of the usual Euclidean algorithm (and harder to analyze as well). 
In terms of continued fractions, cost is the sum
of the partial quotients instead of their number.  A heuristic argument 
suggests a median word length of $O(\log N \log\log N)$
for a matrix with entries in $[-N,N]$.  
By contrast, by a result of D.~Knuth and A.~Yao \KY,
the {\it mean} word length is $O(\log^2 N)$.
The difference between median and mean is due to the fact that a few 
matrices require very long words.  In particular, 
the word length is guaranteed to be large
if the largest matrix entry is much larger in absolute value 
than
the smallest.  The obvious ways of lifting to $\SL_2(\Z)$ nearly always 
produce such unbalanced matrices.  
For example, if we lift the entries of the first row to elements of 
$[0,p-1]$
 and then lift the remaining entries to integers of minimal absolute value,
they will typically be of order $O(p^2)$.  

To avoid this difficulty, we turn the problem around and ask for elements of 
$\SL_2(\Z)$ which {\it can} be represented by short words in our generators.
Let $a$ and $d$ denote integers between $1$ and $p$ with mutually inverse
reductions (mod $p$).
Set $c=p$, so $b=(ad-1)/p$.  For most choices of 
$a$, $a/p$ has a continued fraction expansion with partial quotient sum 
$O(\log p\log\log p)$.
To show this, one must justify the heuristic estimate mentioned above for the sum of the 
partial quotients of a random fraction of fixed demominator $p$.  We do 
this by an elementary argument suggested by the circle 
method.  

In the above construction, $b$ and $d$ are determined by $a$.  To 
eliminate the
dependence on $b$, we use the identity 

$$\mx ab0d \mx1101 \mx ab0d^{-1}=\mx1 {a^2}01.$$
This provides a large number of unitriangular matrices with word 
representations of length $O (\log p \log\log p)$, from which one can 
easily construct all 
elements of $\SL_2(\F_p)$.  

It may be worth noting that the analogous problem for ${\rm SU}(2)$ has recently been solved:
given a fixed finite set of topological generators, to approximate
a given $\alpha\in{\rm SU}(2)$ with error $\epsilon$ by a word of polylog length
in polylog time.  A solution using iterated commutators was discovered
independently by R.~Solovay and A.~Kitaev \CN~App.~3.

I would like to acknowledge the hospitality of the Hebrew University where 
this work was done.   Peter Sarnak first called my attention to the 
problem of efficiently constructing short word representations
for $\SL_2(\F_p)$.  He also made a number of helpful comments 
on an earlier version of this paper.  I enjoyed a number of
stimulating conversations with Alex Lubotzky on this problem.
It gives me great pleasure to thank them both. 

\bigskip
We begin with a careful analysis of the performance of the subtractive 
Euclidean algorithm.  For terminology,
notation, and basic facts related to continued fraction expansions, 
we refer to \HW\ and \K.

\proc{Definition}{An element $\mx abcd $ of $\SL_2 (\Z)$ is {\emph 
left-dominated} if $a, b, c, d\ge 0$ and $a+c\ge b+d$.  }

\proc{Lemma}{If $\mx abcd$ is left-dominated, then $a \ge b$.  Moreover, 
unless it is the identity matrix, $c \ge d$.}

\proof If $a < b $, then $d<c$, which is impossible, since the matrix 
entries are non-negative and the determinant is $1$. Similarly, if $c<d$, 
then $b<a$, so $bc\le ad-1$, with equality if and only if $b=c=0$.\qed

\proc{Lemma}{If $\mx abcd$ is a left-dominated matrix other than the 
identity, then if $a\le c$, the matrix $\mx ab{c-a}{d-b}$ is left-dominated; 
otherwise $\mx{a-c}{b-d}cd$ is left-dominated. }

\proof All that remains to be shown is that the entries of the specified 
matrix are non-negative. If $a \le c$, unimodularity implies $d > b$. If 
$a > c$ and $d > b$, unimodularity implies that $\mx abcd $ is the 
identity. \qed

This lemma shows that the elementary row operations needed to reduce a 
left-dominated matrix to the identity can be chosen without reference to 
the right column. It therefore motivates the definition of a function 
$S:\N\times\N\to\N$ as follows: 
$$ S(a,c)=\cases{0 & if $c=0$, \cr
                 S(a-c, c)+1 & if $a > c > 0$, \cr
                 S(a, c-a)+1 & if $c \ge a > 0$.} $$

We have immediately from this definition the following lemma:
\proc{Lemma}{Any left-dominated matrix $\mx abcd$ can be written as a word 
of length $S(a,c)$ in the letters $\mx 1101$ and $\mx 1011$.}

Every positive rational number has exactly two continued fraction 
expansions: 
$${a\over c}=[k_0, k_1, \ldots, k_n]=[k_0, k_1, \ldots, k_n-1, 1].$$
Therefore, we may define $T(a/c)$ to be the sum of the partial quotients 
appearing in a continued fraction expansion of $a/c$.  

\proc{Lemma}{If $a$ and $c$ are relatively prime
positive integers, $T(a/c)=S(a,c)$.}
\proof Immediate by induction.\qed

\medskip
Our object will be to show that for any fixed prime $p$ there exist many 
positive integers $a<p$ such that $S(a,p) $ is not much greater than $\log 
p\log\log p$.  To do this, it will be convenient to break up $T(a/p)$ into 
pieces corresponding to individual partial quotients.  We therefore define 
$T_d(a/p)$ to be equal to $k_i$ if the denominator of $[k_0, k_1, \ldots, 
k_{i-1}]$ is $d$ for some $i\le n$; to be $1$ if $d$ is the denominator of
$[k_0, k_1, \ldots, k_n-1]$; and otherwise to be $0$.   Thus $T_d(a/p)>0$
if and only if some fraction with denominator $d$ is a convergent of $a/p$.
Moreover,
$$\sum_{d=1}^\infty T_d(a/p)=\sum_{d=1}^{p-1} T_d(a/p)=T(a/p)+1. $$ 

The key proposition is as follows:

\proc{Proposition}{For all $\epsilon>0$, there exists a constant $C$ such 
that for all primes $p$, 
$${|\{a\in[1, p-1]\cap\N:\;S(a,p)<C\log p\log\log p\}|\over p-1}
\ge 1-\epsilon. $$}
\proof Let $X=[1, p-1]\cap\N$. 
For $\delta\in (0,1) $ and $d\in \N$ we define the ``major arc'' $Y_d(\delta)$ 
(really a union of major arcs) to be 
the subset of $X$ consisting of $a$ such that 
$$\inf_{b\in\Z}\left|{a\over p}-{b\over d}\right|<{\delta\over d^2}.$$
For $d \ge p$, $Y_d(\delta)$ is empty.  Every $d<p$ 
is relatively prime to $p$, so there is at most one way to represent a 
given integer as $ad-bp$, $a\in X$, $b\in\Z$, (and no way to represent $0$).  
Thus,
$$|Y_d(\delta)|\le {2p\delta\over d}.$$  
We define 
$$ X(\delta)=X\setminus\bigcup_{d=1}^\infty Y_d(\delta),$$     
so 
$$|X(\delta)|\ge p-1-p\left(2\delta\sum_{d=1}^{p-1} 1/d\right)
\ge -1+p(1-2\delta(\log p + 1))\ge -1+p(1-6\delta\log p).$$
On the other hand, the ``minor arc'' contribution satisfies
$$\sum_{a\in X(\delta)}T(a/p)<\sum_{a\in X(\delta)}\sum_{d=1}^{p-1}T_d(a/p).$$
If $p_i/q_i$ denotes the $i$th convergent of $a/p$, then
$$\left|{p_i\over q_i}-{a\over p}\right|< {1\over k_{i+1}q_i^2}.$$
In particular, unless $|bp- ad|< p/d$, $b/d$ cannot be a convergent of 
$a/p$.  Thus, 
$$\sum_{a\in X(\delta)}{T_d(a/p)}\le 2\sum_{e\in [p\delta /d, 
p/d]\cap\N}{p\over de}\le {2\over\delta}-{2p\log\delta \over d}.$$
Summing over $d$, 
$$\sum_{a\in X(\delta)}T(a/p)<{(2p-2)\over\delta} - 2\log\delta(\log p+1)p.  $$
Setting $\delta = {\epsilon\over 12\log p}$, we get $|X(\delta)|\ge 
(1-\epsilon/2)p-1$ and 
$$\sum_{a\in X(\delta)}T(a/p)\le 2p\log p\log\log p+o(p\log p\log\log p).$$

Choosing $C$ sufficiently large, the number of elements $a$ in $X(\delta)$ 
with 
$$T(a/p)>C\log p\log\log p$$
is less than $\epsilon p/2$.  Thus, the 
number of elements in $X$ with $T(a/p)>C\log p\log\log p$ is at most 
$\epsilon p$.\qed

We can now prove the main theorem:

\proof Setting $\epsilon$ to be any constant less than $1/16$, we 
define $C$ as above.  As the number of points on a nonsingular affine 
conic over $\F_p$ is at least $p-1$ and at most $p+1$, for any $y\in\F_p$, 
there are $p/4+O(1)$ representations of $y$ as a sum of quadratic 
residues $x_1$ and $x_2$.  We write $a_i$ for the representative of $\sqrt{x_i}$ 
in $[1, p/2]$.  The number of choices of $x_1$ for which 
$$\sup_i T(a_i/p)>C\log p\log\log p$$ 
is at most $4\epsilon p$, so that if an element $x_1$ of $\F_p$ is chosen at
random,
the probability is at least $1/4-2\epsilon>0$ that $x_1$ is a 
square, $S(a_1, p)\le C\log p\log\log p$, and the same things 
are true for $x_2=y-x_1$ and the unique integer $a_2\in [1, p/2]$ such that 
$a_2^2$ reduces to $x_2$.  Define $d_i$ to be the integer in $[1, p-1]$ 
which reduces to the inverse of the reduction of $a_i$, and set $b_i=(a_i 
d_i-1)/p$. Thus, 
$$\mx {a_1}{b_1}0{d_1}\mx 1101\mx {a_1}{b_1}0{d_1}^{-1}\mx 
{a_2}{b_2}0{d_2}\mx 1101\mx {a_2}{b_2}0{d_2}^{-1}$$
can be written as a word of length at most $4C\log p\log\log p+2$, and its 
(mod $p$) reduction is $\mx 1y01$.  

For square roots,  we use Shanks's algorithm (\Shanks, \Lind),
which is probabilistic and polylogarithm.  Note that
one has a deterministic square root algorithm when $p\equiv 3\ 
\hbox{(mod $4$)}$, but nevertheless, our algorithm remains nondeterministic
since it depends on how many tries are needed before we find a good $x_1$. 

Applying transpose, we can likewise find words of length $O(\log p\log\log p)$
for lower unitriangular matrices.  Since every matrix which is not 
upper triangular can be written 
$$\mx 1{y_1}01\mx 10{y_2}1\mx 1{y_3}01, $$
every matrix in $\SL_2(\F_p)$ can be written as a product of at most four 
upper or lower unitriangular matrices.  Therefore, for every constant $c_4 < (1/4)^4$, we can find 
$c_1$, $c_2$, and $c_3$ satisfying the conditions of the theorem.  To deal with $c_4\ge 4^{-4}$,
we use repeated independent trials of the algorithm.

\qed
\biblio

\end